\documentclass[12pt,amssymb]{amsart}
\usepackage{amssymb}
\usepackage[mathscr]{eucal}

\newcommand{\Z}{{\mathbb Z}}
\newcommand{\C}{{\mathbb C}}
\newcommand{\Q}{{\mathbb Q}}
\newcommand{\cL}{\mathcal{L}}
\newcommand{\R}{{\mathbb R}}

\newcommand{\Tr}{\mathrm{Tr}\: }
\newcommand{\cK}{\mathcal{K}}

\newcommand{\cR}{\mathcal{R}}
\newcommand{\cE}{\mathcal{E}}

\newcommand{\md}{{\rm mod}\ }
\newcommand{\fX}{\mathfrak{X}}
\newcommand{\cO}{\mathcal{O}}

\newcommand{\Ga}{\mathrm{Gal}}
\newcommand{\fH}{\mathfrak{H}}

\newcommand{\cG}{\mathcal{G}}

\begin{document}

\title[Number-theoretic techniques]{Number-theoretic techniques in the theory of Lie groups and
differential geometry}

\author[Prasad]{Gopal Prasad}

\author[Rapinchuk]{Andrei S. Rapinchuk}

\address{Department of Mathematics, University of Michigan, Ann
Arbor, MI 48109}

\email{gprasad@umich.edu}

\address{Department of Mathematics, University of Virginia,
Charlottesville, VA 22904}

\email{asr3x@virginia.edu}

\maketitle

The aim of this article is to give a brief survey of the results
obtained in the series of papers \cite{PR1}--\cite{PR5}. These papers
deal with a variety of problems, but have a common feature:\:they all
rely in a very essential way on number-theoretic
techniques (including $p$-adic techniques), and use  results from
algebraic and transcendental number theory. The fact that 
number-theoretic techniques
turned out to be crucial for tackling certain problems originating in the theory of
(real) Lie groups and differential geometry was very exciting. We hope that these 
techniques will become an integral part of the repertoire of mathematicians working in these areas.

To keep the size of this article within a reasonable limit, we will
focus primarily on the paper \cite{PR5}, and briefly mention the results of \cite{PR1}--\cite{PR4} 
and some other related results in
the last section. The work in \cite{PR5}, which was originally motivated by
questions in differential geometry dealing with {\it length-commensurable} 
and {\it isospectral} locally symmetric spaces (cf.\,\S 1), led us to define a new
relationship between Zariski-dense subgroups of a simple (or
semi-simple) algebraic group which we call {\it weak
commensurability} (cf.\,\S 2). The results of \cite{PR5} give an
almost complete characterization of weakly commensurable arithmetic
groups, but there remain quite a few natural questions (some of
which are mentined below) for general Zariski-dense subgroups. We
hope that the notion of weak commensurability will be
useful in investigation of (discrete) subgroups of Lie groups,
geometry and ergodic theory.

\section{Length-commensurable and isospectral manifolds}

Let $M$ be a Riemannian manifold. In differential geometry, one
associates to $M$ the following sets of data: the {\it length
spectrum} $\cL(M)$ (the set of lengths of all closed geodesics {\it
with} multiplicities), the {\it weak length spectrum} $L(M)$ (the
set of lengths of all closed geodesics {\it without}
multiplicities), the {\it spectrum of the Laplace operator} $\cE(M)$
(the set of eigenvalues of the Laplacian $\Delta_M$ with
multiplicities). The fundamental question is to what extent do
$L(M),$ $\cL(M)$ and $\cE(M)$ determine $M?$ In analyzing this
question, the following terminology will be used: two Riemannian
manifolds $M_1$ and $M_2$ are said to be {\it isospectral} if
$\cE(M_1) = \cE(M_2),$ and {\it iso-length-spectral} if $\cL(M_1) =
\cL(M_2).$

First, it should be pointed out that the conditions like isospectrality,
iso-length-spectrality are related to each other. For example, for
compact hyperbolic 2-manifolds $M_1$ and $M_2,$ we have $\cL(M_1) =
\cL(M_2)$ if and only if $\cE(M_1) = \cE(M_2)$ (cf.\,\cite{Mc}), and
two hyperbolic 3-manifolds are isospectral if and only if they have
the same {\it complex-length} spectrum (for its definition see the footnote 
later in this section),
cf.\,\cite{BB} or \cite{Ga}. Furthermore, for compact locally
symmetric spaces $M_1$ and $M_2$ of nonpositive curvature, if
$\mathcal{E}(M_1) = \mathcal{E}(M_2)$, then $L(M_1) = L(M_2)$ (see
\cite{PR5}, Theorem 10.1). (Notice that all these results rely on
some kind of trace formula.)

Second, neither of $\mathcal{L}(M),$ $L(M)$ or $\mathcal{E}(M)$
determines $M$ up to isometry. In fact, in 1980, Vign\'eras \cite{Vi}
constructed examples of isospectral, but nonisometric, hyperbolic 2
and 3-manifolds. This construction relied on arithmetic properties
of orders in a quaternion algebra $D$. More precisely, her crucial
observation was that it is possible to choose $D$ so that it contains
orders $\cO_1$ and $\cO_2$ with the property that the corresponding
groups $\cO_1^{(1)}$ and $\cO_2^{(1)}$ of elements with reduced norm
one are not conjugate, but their closures in the completions are
conjugate, for all nonarchimedean places of the center. Five years
later, Sunada \cite{Su} gave a very general, and purely
group-theoretic, method for constructing isospectral, but
nonisometric, manifolds. His construction goes as follows: {\it
\:Let $M$ be a Riemannian manifold with the fundamental group
$\Gamma := \pi_1(M).$ Assume that $\Gamma$ has a finite quotient $G$
with the following property: there are subgroups $H_1 , H_2$ of $G$
such that $\vert C \cap H_1 \vert = \vert C \cap H_2 \vert$ for all
conjugacy classes $C$ of $G.$ Let $M_i$ be the finite-sheeted cover
of $M$ corresponding to the pull-back of $H_i$ in $\Gamma.$ Then
(under appropriate assumptions), $M_1$ and $M_2$ are nonisometric
isospectral (or iso-length-spectral) manifolds.}

Since its inception, Sunada's method and its variants have been
used to construct examples of nonisometric
manifolds with same invariants. In particular, Alan Reid
\cite{Re} constructed examples of nonisometric iso-length-spectral
hyperbolic 3-manifolds, and last year, in a joint paper \cite{LMNR},
Leninger, McReynolds, Neumann and Reid gave examples of hyperbolic
manifolds with the same weak length spectrum, but different volumes.
These, and other examples, demonstrate that it is not possible 
to characterize  Riemannian manifolds (even hyperbolic ones) up to
isometry by their spectrum or length spectrum.
On the other hand, it is worth noting that the manifolds furnished
by Vign\'eras, and the ones obtained using Sunada's method are always {\it
commensurable,} i.e., have a common finite-sheeted covering. This
suggests that the following is perhaps a more reasonable question.

\vskip2mm

\noindent {\bf Question 1:} {\it Let $M_1$ and $M_2$ be two
(hyperbolic) manifolds (of finite volume or even compact). Suppose
$L(M_1) = L(M_2).$ Are $M_1$ and $M_2$ necessarily commensurable?}

\vskip2mm

\noindent (Of course, the same question can be asked for other
classes of manifolds, e.g.\:for general locally symmetric spaces of finite volume.)

\vskip2mm

The answer even to this modified question turns out to be ``no" in
general: Lubotzky, Samuels and Vishne \cite{LSV} have given examples
of  isospectral (hence, with same weak length spectrum) compact locally
symmetric spaces that are {\it not} commensurable. At the same time,
some positive results have emerged. Namely,
Reid \cite{Re} and Chinburg, Hamilton, Long and Reid \cite{CHLR}
have given a positive answer to Question 1 for arithmetically defined
hyperbolic 2- and 3-manifolds, respectively. Our results in
\cite{PR5} provide an almost complete answer to Question 1 for
arithmetically defined locally symmetric spaces of arbitrary absolutely simple
Lie groups. In fact, in \cite{PR5} we analyze when two
locally symmetric spaces are commensurable given that they satisfy a
much weaker condition than iso-length-spectrality, which we termed
{\it length-commensurability.} We observe that not only does the use
of this condition produce stronger results, but the condition itself
is more suitable for analyzing Question 1 as it allows one to
replace the manifolds under consideration with commensurable
manifolds.

\vskip2mm

\noindent {\bf Definition.} Two Riemannian manifolds $M_1$ and $M_2$ are said to be {\it
length-commensurable} if $\Q \cdot L(M_1) = \Q \cdot L(M_2).$

\vskip2mm

Now, we are in a position to formulate precisely the question which
is central to \cite{PR5}.

\vskip2mm

\noindent {\bf Question 2:} {\it Suppose $M_1$ and $M_2$ are
length-commensurable. Are they commensurable?}

\vskip2mm

In \cite{PR5}, we have been able to answer this question for 
arithmetically defined locally symmetric spaces of absolutely simple Lie groups. The precise
formulations will be given in \S 3, after introducing appropriate definitions.
The following theorem, however, is fully representative of these
results.

\vskip3mm

\noindent {\bf Theorem.} (1) {\it Let $M_1$ and $M_2$ be two {\rm
arithmetically defined hyperbolic} manifolds of {\rm even}
dimension. If $M_1$ and $M_2$ are not commensurable, then, after a
possible interchange of $M_1$ and $M_2,$ there exists $\lambda_1 \in
L(M_1)$ such that for {\rm any} $\lambda_2 \in L(M_2),$ the ratio
$\lambda_1/\lambda_2$ is {\rm transcendental}. In
particular, $M_1$ and $M_2$ are not length-commensurable.}

\vskip2mm

\noindent (2) {\it For any dimension} $d \equiv 1(\md 4),$ {\it
there exist length-commensurable, but not commensurable,
arithmetically defined hyperbolic $d$-manifolds.}

\vskip2mm

We have proved similar results for arithmetically defined locally
symmetric spaces of absolutely simple real Lie groups of all types; see \cite{PR5} and \cite{PR6}. For example, for hyperbolic spaces modeled on Hamiltonian
quaternions we have an assertion similar to (1) (i.e., Question 2
has an affirmative answer); but for complex hyperbolic spaces we
have an assertion similar to (2) (i.e., Question 2 has a negative
answer).

\vskip1mm

The key ingredient of our approach is the new notion of {\it weak
commensurability} of Zariski-dense subgroups of an algebraic group,
and the {relationship between the length-commensurability} of
locally symmetric spaces and the {weak commensurability} of their
fundamental groups. To motivate the definition of weak
commensurability, we consider the following simple example.

\vskip1mm

Let $\fH = \{ x + iy \: \vert\: y > 0 \}$ be the upper half-plane
with the standard hyperbolic metric $ds^2 = y^{-2}(dx^2 + dy^2).$
Then $t\mapsto ie^t$ is a geodesic in $\fH,$ whose piece $\tilde{c}$
connecting $i$ to $ai,$ where $a > 1$, has length $\ell(\tilde{c}) =
\log a.$ Now, let $\Gamma \subset SL_2(\R)$ be a discrete
torsion-free subgroup, and $\pi \colon \fH \to \fH/\Gamma$ be the
canonical projection. If $c:=\pi(\tilde{c})$ is a closed geodesic 
in $\fH/\Gamma$ (traced once), then it is not difficult to see that
for $\lambda = \sqrt{a},$ the element
$$\gamma =
\left(\begin{array}{cc} \lambda& 0 \\ 0 & \lambda^{-1} \end{array}
\right)$$ lies in $\Gamma$. Then the length $\ell(c)$ of $c$ equals 
$\log a = 2\log \lambda.$  This shows that the lengths of closed geodesics in the hyperbolic
2-manifold $\fH/\Gamma$ are (multiples of) the logarithms of the
eigenvalues of semi-simple elements of the fundamental group
$\Gamma$ (cf.\,\S 3 below,  and Proposition 8.2 in \cite{PR5} for a
general statement that applies to arbitrary locally symmetric
spaces)\footnote{In the above construction if we replace ${\rm SL}_2(\R)$ with
${\rm SL}_2(\C)$, then the collection of (principal values) of the
logarithms of the eigenvalues of semi-simple elements is known as
the {\it complex length spectrum} of the corresponding hyperbolic
3-manifold.}. Furthermore, let $c_i$ for $i = 1 , 2,$ be closed geodesics
in $\fH/\Gamma$ that in the above notation correspond to
semi-simple elements $\gamma_i \in \Gamma$ having the eigenvalue 
$\lambda_i > 1.$ Then
$$\ell(c_1)/\ell(c_2) = m/n \ \ \Leftrightarrow \ \ \lambda_1^n = \lambda_2^m.$$
Notice that the condition on the right-hand side can be reformulated 
as follows: {\it If $T_i$ is a torus of $\mathrm{SL}_2$ such that $\gamma_i\in T_i(\R)$, then there exist $\chi_i \in X(T_i)$ with 
$$\chi_1(\gamma_1) = \chi_2(\gamma_2) \neq 1.$$}

The above discussion suggests the following.

\vskip1mm

\noindent {\bf Definition.} Let $G$ be a semi-simple algebraic group
defined over a field $F.$ Two semi-simple elements $\gamma_1 , \gamma_2 \in G(F)$
are {\it weakly commensurable} if, for $i=1,\,2$, there exist maximal $F$-tori $T_i$,
and characters $\chi_i \in X(T_i)$,  such that $\gamma_i\in T_i(F)$, and $$\chi_1(\gamma_1) = \chi_2(\gamma_2) \neq 1.$$

\vskip2mm

As we have seen, weak commensurability adequately reflects
length-commensurability of hyperbolic 2-manifolds. In fact, it
remains relevant for length-commensurability of arbitrary locally
symmetric spaces. This is easy to see for rank one spaces but is
less obvious for higher rank spaces - cf.\,\S 3 below,  and \cite{PR5}, \S 8.

\section{Weakly commensurable arithmetic subgroups}

We observe that for $G \neq \mathrm{SL}_2,$ weak commensurability of
$\gamma_1 , \gamma_2 \in G(F)$ may not relate these elements to each other in a
significant way (in particular, $F$-tori $T_i$ of $G$ containing
these elements may be very different). So, to get meaningful
consequences of weak commensurability, one needs to extend this
notion from individual elements to ``large"  (in particular,
Zariski-dense) subgroups.

\vskip2mm

\noindent {\bf Definition.} Two (Zariski-dense) subgroups $\Gamma_1
, \Gamma_2$ of $G(F)$ are {\it weakly commensurable} if every
semi-simple element $\gamma_1 \in \Gamma_1$ of infinite order is
weakly commensurable to some semi-simple element $\gamma_2 \in
\Gamma_2$ of infinite order, and vice versa. \vskip2mm

\vskip2mm

It was discovered in \cite{PR5} that weak commensurability has some
important consequences even for completely general finitely generated 
Zariski-dense subgroups. For
simplicity, we will assume henceforth that all our fields are of
characteristic zero. To formulate our first result, we need one
additional notation: given a subgroup $\Gamma$ of  $G(F),$ where
$G$ is an absolutely simple algebraic $F$-group, we let $K_{\Gamma}$ denote the
subfield of $F$ generated by the traces $\mathrm{Tr}\: \mathrm{Ad}\:
\gamma$ for all $\gamma \in \Gamma$, where $\mathrm{Ad}$ denotes the
adjoint representation of $G.$ We recall that according to a result
of Vinberg \cite{Vn}, for a Zariski-dense subgroup $\Gamma$ of $G$, 
the field $K_{\Gamma}$ is precisely the field of definition of
$\mathrm{Ad}\: \Gamma,$ i.e.,\, it is  the minimal subfield $K$ of $F$ 
such that all elements of $\mathrm{Ad}\: \Gamma$ can be
represented simultaneously by matrices with entries in $K,$ in a
certain basis of the Lie algebra $\mathfrak{g}$ of $G.$

\vskip2mm

\noindent {\bf Theorem A.} {\it Let $\Gamma_1$ and $\Gamma_2$ be two
finitely generated Zariski-dense subgroups of $G(F).$ If $\Gamma_1$
and $\Gamma_2$ are weakly commensurable, then $K_{\Gamma_1} =
K_{\Gamma_2}.$}

\vskip2mm

Much stronger results are available for the case of arithmetic
subgroups. To formulate these, we need to describe  the terminology we
use regarding arithmetic subgroups. Let $G$ be a semi-simple
algebraic group over a field $F$ of characteristics zero. Suppose we
are given:

\vskip2mm

 $\bullet$ {a number field $K$ contained in $F$;}

\vskip1mm

 $\bullet$ {a (finite) subset $S$ of places of $K$ containing all the archimedean places;}

\vskip1mm

 $\bullet$ {a $K$-form $G_0$ of $G,$ i.e., a group $G_0$ defined over $K$ such that $_FG_0 \stackrel{\iota}{\simeq} G$ over
$F.$}

\vskip2mm

\noindent Then subgroups of $G(F)$ commensurable with the image of
the natural embedding $G_0(\cO_K(S))\hookrightarrow G(F)$ induced by
$\iota,$  where $\cO_K(S)$ is the ring of $S$-integers in $K$, are
by definition $(G_0 , K , S)$-{\it arithmetic} subgroups. Notice
that in this definition we {\it do} fix an embedding of $K$ into $F$
(in other words, isomorphic, but distinct, subfields of $F$ are
treated as different fields), but we {\it do not} fix an
$F$-isomorphism $\iota,$ so by varying it we generate a class of
subgroups {\it invariant under $F$-automorphisms.} For this reason,
by ``commensurability" we will mean ``commensurability {\it up to
$F$-isomorphism,}" i.e., two subgroups $\Gamma_1$ and $\Gamma_2$ of  $G(F)$ will
be called commensurable if there exists an $F$-automorphism
$\varphi$ of $G$ such that $\varphi(\Gamma_1)$ and $\Gamma_2$ are
commensurable in the usual sense, viz. their intersection has finite
index in both of them. Another convention is that $S$ will be
assumed to contain no nonarchimedean places $v$ such that $G_0$ is
$K_v$-anisotropic (this assumption enables us to recover $S$
uniquely from a given $S$-arithmetic subgroup).

\vskip2mm

{\it The group $G$ in Theorems B--F is assumed to be absolutely simple.}

\vskip2mm

\noindent {\bf Theorem B.} {\it Let $\Gamma_i$ be a
Zariski-dense $(G_i , K_i , S_i)$-arithmetic subgroup of $G(F)$ for $i = 1 , 2.$
If $\Gamma_1$ and $\Gamma_2$ are weakly commensurable, then $K_1 =
K_2$ and $S_1 = S_2.$}

\vskip2mm

One shows that $\Gamma_1$ and $\Gamma_2$ as in Theorem B are
commensurable if and only if $K_1 = K_2,$ $S_1 = S_2$ and $G_1
\simeq G_2$ over $K := K_1 = K_2$ (cf.\,Proposition 2.5 in
\cite{PR5}). So, according to Theorem B, the weak commensurability
of $\Gamma_1$ and $\Gamma_2$ implies that the first two of these
three conditions do hold. In general, however, $G_1$ and $G_2$
do not have to be $K$-isomorphic. Our next theorem describes the
situations where it can be inferred that $G_1$ and $G_2$ are $K$-isomorphic.

\vskip2mm

\noindent {\bf Theorem C.} {\it Suppose $G$ is not of type $A_n$ $(n
> 1),$ $D_{2n+1}$ $(n \geqslant 1)$ or $E_6.$ If $G(F)$ contains
Zariski-dense weakly commensurable $(G_i , K , S)$-arithmetic
subgroups $\Gamma_i$ for $i = 1 , 2,$ then $G_1 \simeq G_2$ over
$K,$ and hence $\Gamma_1$ and $\Gamma_2$ are commensurable up to an
$F$-automorphism of $G.$}

\vskip2mm

In the general case, we have the following finiteness result.

\vskip2mm

\noindent {\bf Theorem D.} {\it Let $\Gamma_1$ be a Zariski-dense
$(G_1 , K , S)$-arithmetic subgroup of $G(F).$ Then the set of
$K$-isomorphism classes of $K$-forms $G_2$ of $G$ such that $G(F)$
contains a Zariski-dense $(G_2 , K , S)$-arithmetic subgroup weakly
commensurable to $\Gamma_1$ is finite. In other words, the set of all  $(K , S)$-arithmetic
subgroups of $G(F)$ which are weakly commensurable to a given $(K ,
S)$-arithmetic subgroup is a union of finitely many commensurability
classes.}

\vskip2mm

Note that for the types $A_n$ ($n > 1$), $D_{2n+1}$ and $E_6$ excluded in
Theorem C, the number of commensurability classes in Theorem D may
not be bounded by an absolute constant depending, say, on $G,$ $K$ and
$S:$ as one varies $\Gamma_1$ (or, equivalently, $G_1$), this number
changes and typically grows to infinity. To explain what happens for
groups of these types, let us consider the following example.

\vskip2mm

Fix any $n > 1$ and pick four nonarchimedean places
$v_1,\,v_2,\,v_3,\, v_4 \in V^K.$ Next, consider central division
$K$-algebras $D_1$ and $D_2$ of degree $d = n + 1 > 2$ with local
invariants $(\in \Q/\Z):$
$$
n_v^{(1)} = \left\{\begin{array}{rcl} 0 &, & v \neq v_i, i\leqslant 4\\
1/d &, & v = v_1 \ \text{or} \  v_2 \\ -1/d &, & v = v_3 \ \text{or}
\ v_4
\end{array} \right.
$$
and
$$
n_v^{(2)} = \left\{\begin{array}{rcl} 0 &, & v \neq v_i, i\leqslant 4\\
1/d &, & v = v_1 \ \text{or} \ v_3 \\ -1/d &, & v = v_2 \ \text{or}
\ v_4.
\end{array} \right.
$$
\vskip1mm

\noindent Then the algebras $D_1$ and $D_2$ are neither isomorphic nor
anti-isomorphic, implying that the algebraic groups $G_1 =
\mathrm{SL}_{1 , D_1}$ and $G_2 = \mathrm{SL}_{1 , D_2}$ (which are
anisotropic inner forms of type $A_n$) are not $K$-isomorphic. On
the other hand, $D_1$ and $D_2$ have exactly the same maximal
subfields, which means that $G_1$ and $G_2$ have the same maximal
$K$-tori. It follows that for any $S,$ the corresponding
$S$-arithmetic subgroups are weakly commensurable, but not
commensurable. Furthermore, by increasing the number of places in
this construction, one can construct an arbitrarily large number of
central division $K$-algebras of degree $d$ with the above
properties. Then the associated $S$-arithmetic groups will all be
weakly commensurable, but will constitute an arbitrarily large
number of commensurability classes.

In \cite{PR5}, Example 6.6, we described how a similar construction
can be given  for {\it some} outer form of type $A_n$ (i.e.,\, for
special unitary groups), at least when $d = n + 1$ is odd. The
restriction on $d$ is due to the fact that our argument
relies on  a local-global principle for embedding fields with an
involutive automorphisms into an algebra with an involution of the
second kind (Proposition A.2 in \cite{PR0}), which involves some
additional assumptions. Recently, we have been able to remove any
restrictions in the local-global principle (unpublished), so the
construction can in fact be implemented for all $d.$

However, no construction of nonisomorphic $K$-groups with the same
$K$-tori was known for types $D_n$ and $E_6.$ We have given a
construction, using Galois cohomology, which works uniformly for
types $A_n,$ $D_{2n+1}$ and $E_6$ (cf. \cite{PR5}, \S 9). Towards
this end, we established a new local-global principle for the
existence of an embedding of a given $K$-torus as a maximal torus in
a given absolutely simple simply connected $K$-group. This
construction, of course, allows one to produce examples of
noncommensurable weakly commensurable  $S$-arithmetic subgroups in
groups of types $A_n,$ $D_{2n+1}$ and $E_6,$ and in fact, show that
the number of commensurability classes is unbounded. This construction may also be
useful elsewhere, for example, in the Langlands
program.

\vskip2mm

Even though the definition of weak commensurability
involves only semi-simple elements, it detects the presence of
unipotent elements; in fact it detects $K$-rank.

\vskip2mm

\noindent {\bf Theorem E.} {\it Assume that $G(F)$ contains
Zariski-dense weakly commensurable $(G_1 , K , S)$- and $(G_2 , K ,
S)$-arithmetic subgroups. Then the Tits indices of $G_1/K$ and
$G_2/K$ are isomorphic. In particular, $\mathrm{rk}_K\: G_1 =
\mathrm{rk}_K\: G_2.$}

\vskip2mm

The above results provide an almost complete picture of weak
commensurability among $S$-arithmetic subgroups. In view of
the connection of weak commensurability with length-commensurability of
locally symmetric spaces (cf.\,\S 3), one would
like to extend these results to not necessarily
arithmetic Zariski-dense subgroups. We conclude this section with an arithmeticity
theorem in which only one subgroup is assumed to be arithmetic, and
a discussion of some open questions.
\vskip2mm

\noindent{\bf Theorem F.} {\it Let $G$ be an absolutely simple
algebraic group over  a nondiscrete locally compact field $F$. Let
$\Gamma_1$ and $\Gamma_2$ be two lattices in $G(F),$ with $\Gamma_1$
$(K,S)$-arithmetic. If $\Gamma_1$ and $\Gamma_2$ are weakly
commensurable, then $\Gamma_2$ is also $(K,S)$-arithmetic.}

\vskip2mm

\noindent {\bf Remarks:} (i) The assumption that {\it both}
$\Gamma_1$ and $\Gamma_2$ be lattices cannot be omitted. For
example, let $\Gamma \subset SL_2(\Z)$ be a torsion-free subgroup of
finite index, and $\Gamma^n$ be the subgroup generated by the $n$-th
powers of elements in $\Gamma$. Then  $\Gamma^n$ is weakly
commensurable with $\Gamma$ for all $ n.$ On the other hand,
$[\Gamma : \Gamma^n] = \infty$ for all sufficiently large $n,$ and
then $\Gamma^n$ is not arithmetic. The same remark applies to 
all hyperbolic groups. However, we do not know what happens in the higher rank situation.

\vskip2mm

\noindent (ii) The case of lattices in products of real and $p$-adic
groups has not been fully investigated.

\vskip2mm

\noindent (iii) Yet another interesting open question is whether or
not the discreteness of one of the two weakly commensurable
subgroups $\Gamma_1 , \Gamma_2$ of  $G(F)$ implies the
discreteness of the other (here $F$ is a locally compact nondiscrete
field).

\vskip2mm

Further analysis of weak commensurability of general Zariski-dense
subgroups of $G(F)$ for an arbitrary field $F$ would require
information about classification of forms of $G$ over general
fields, which is not yet available. For example, even the
following basic question seems to be open.

\vskip1.5mm

\noindent {\bf Question 3:} {\it Let $D_1$ and $D_2$ be two quaternion
division algebras over a finitely generated field $K.$ Assume that
$D_1$ and $D_2$ have the {\rm same maximal subfields}. Are they
isomorphic?}

\vskip1.5mm

M.\,Rost has informed us that over large fields (like those used in the
proof of the Merkurjev-Suslin theorem) the answer can be ``no"
(apparently, the same observation was independently made by
A.~Wadsworth and some other people). But for finitely generated
fields (note that the fields arising in the investigation of weakly
commensurable finitely generated subgroups are finitely generated),
the answer is unknown. Furthermore, if the answer turns out to
be ``no", we would like to know if the number of isomorphism
classes of quaternion algebras over a given finitely generated field, 
and containing the same maximal subfields is finite (this may be useful for
extending the finiteness result of Theorem D to such nonarithmetic
subgroups as the fundamental groups of general compact Riemann surfaces).

\vskip2mm

\section{Length-commensurable locally symmetric spaces}

\vskip3mm

Let $G$ be a connected semi-simple algebraic $\R$-group, $\cG
= G(\R).$ We let $\cK$ denote a maximal  compact subgroup of $\cG,$
and let  $\mathfrak{X} =\cK \backslash \cG$ be the corresponding
symmetric space of $\cG$. For a discrete torsion-free subgroup
$\Gamma$ of $\cG,$ we let $\fX_{\Gamma} := {\fX}/\Gamma$ denote
the locally symmetric space with the fundamental group $\Gamma.$ We
say that ${\fX}_{\Gamma}$ is {\it arithmetically defined} if
$\Gamma$ is arithmetic (with $S$ the set of archimedean places of $K$) in the sense
specified in \S 2. Notice that $\mathfrak{X}_{\Gamma_1}$ and
$\mathfrak{X}_{\Gamma_2}$ are commensurable as manifolds (i.e.,\,
have a common finite-sheeted cover) if and only if  $\Gamma_1$ and
$\Gamma_2$ are commensurable up to $\R$-automorphism of $G.$

\vskip1mm

Our goal now is to relate length-commensurability of locally
symmetric spaces to weak commensurability of their fundamental
groups. We need to recall some basic facts about closed geodesics in
$\fX_{\Gamma}$ (cf.\,\cite{PR4}, or \cite{PR5}, \S 8). The closed geodesics on $\fX/\Gamma$ correspond to semi-simple elements of $\Gamma$.  For a semi-simple element $\gamma\in \Gamma$, let $c_{\gamma}$ be the closed geodesic corresponding to $\gamma$. Its length is given by the
following formula (see \cite{PR5}, Proposition 8.2).
\begin{equation}\label{E:3-1}
\ell_{\Gamma}(c_{\gamma})^2 = (1/n_{\gamma}^2) \left(\sum (\log
\vert \alpha(\gamma) \vert)^2\right),
\end{equation}
where $n_{\gamma}$ is an integer, and the sum is over all roots
$\alpha$ of $G$ with respect to a maximal $\R$-torus $T$ such that
$\gamma\in T(\R).$ (We notice that for the upper half-plane $\fH = {\rm SO}_2 \backslash 
{\rm SL}_2 (\R)$ this metric differs from the standard
hyperbolic metric, considered in \S 1, by a factor of $\sqrt{2},$
which, of course, does not affect length commensurability.)

For our purposes, we need to recast (\ref{E:3-1}) using the notion
of a {\it positive real character}. Given a real torus $T,$ a real
character $\chi$ of $T$ is called {\it positive} if $\chi(t) > 0$
for all $t \in T(\R).$ We notice that for any character $\chi$ of
$T$ we have
$$\vert \chi(t) \vert^2 = \chi(t) \cdot
\overline{\chi(t)} = (\chi + \overline{\chi})(t) = \chi_0(t),$$
where $\chi_0$ is a positive real character. Hence,
\begin{equation}\label{E:3-2}
\ell_{\Gamma}(c_{\gamma})^2 = (1/n_{\gamma}^2)\sum_{i = 1}^p s_i
(\log \chi^{(i)}(\gamma))^2,
\end{equation}
where $s_i \in \Q$, and $\chi^{(i)}$ are positive real characters.

\vskip2mm

The right-hand side of (\ref{E:3-2}) is easiest to analyze when
$\mathrm{rk}_{\R}\: G = 1,$ which we will now assume. Let $\chi$ be
a generator of the group of positive real characters of a maximal $\R$-torus
$T$ containing $\gamma.$ Then
$$\ell_{\Gamma}(c_{\gamma}) = (s/n_{\gamma}) \cdot \vert
\log \chi(\gamma) \vert,$$ where $s$ is independent of $\gamma$ and
$T$ (because any two maximal $\R$-tori of real rank one are conjugate to each other  by an element
of $G(\R)$). So, if $\gamma_1 \in \Gamma_1$ and $\gamma_2 \in
\Gamma_2$ are {\it not} weakly commensurable, then
\begin{equation}\label{E:3-3}
\ell_{\Gamma_1}(c_{\gamma_1})/\ell_{\Gamma_2}(c_{\gamma_2}) =
(n_{\gamma_2}/n_{\gamma_1}) \cdot \left(\pm \frac{\log
\chi_1(\gamma_1)}{\log \chi_2(\gamma_2)} \right) \notin \Q.
\end{equation}
Therefore, if $\Gamma_1$ and $\Gamma_2$ are {\it not weakly
commensurable}, then $\fX_{\Gamma_1}$  and $\fX_{\Gamma_2}$ are {\it
not length-commensurable.} Thus, the connection noted in \S 1 for
hyperbolic 2-manifolds remains valid for arbitrary locally symmetric
spaces of rank one. In fact, one can make a stronger statement
assuming that $\Gamma_1$ and $\Gamma_2$ are arithmetic (or, more
generally, can be conjugated into ${\rm SL}_n(\overline{\Q})$). Then
$\chi_i(\gamma_i) \in \overline{\Q}^{\times}$ for \ \  $i = 1 , \,
2.$ But according to a theorem proved independently by Gelfond and
Schneider in 1934, if $\alpha$ and $\beta$ are algebraic numbers
such that ${\log \alpha}/{\log \beta}$ is irrational, then it is
transcendental over $\Q$ (cf. \cite{Ba}). So, it follows from
(\ref{E:3-3}) that if $\Gamma_1$ and $\Gamma_2$ are as above, and
$\gamma_1\in\Gamma_1$ and $\gamma_2\in\Gamma_2$ are not weakly commensurable, then
$$\ell_{\Gamma_1}(c_{\gamma_1})/\ell_{\Gamma_2}(c_{\gamma_2})$$ is
transcendental over $\Q.$

\vskip3mm

\noindent To relate length-commensurability of locally symmetric spaces of higher rank 
with the notion of weak commensurability of their fundamental groups, we
need to invoke the Schanuel's Conjecture from transcendental number
theory (cf.\,\cite{Ba}).

\vskip2mm

\noindent {\bf Schanuel's conjecture.} {\it If $z_1, \ldots , z_n
\in \C$ are linearly independent over $\Q$, then the transcendence
degree over $\Q$ of the field generated by $$z_1, \ldots , z_n; \
e^{z_1}, \ldots , e^{z_n}$$ is $\geqslant n.$}

\vskip2mm

What we need is the following corollary of Schanuel's conjecture.
Let $\alpha_1, \ldots , \alpha_n$ be nonzero {\it algebraic}
numbers, and set $z_i = \log \alpha_i.$ Applying Schanuel's
conjecture, we obtain that $\log \alpha_1, \ldots , \log \alpha_n$
are {\it algebraically independent}  as soon as they are {\it
linearly independent} (over $\Q$), i.e.,\, whenever $\alpha_1, \ldots ,
\alpha_n$ are multiplicatively  independent.

\vskip2mm

Before we proceed, we would like to point out that our results for
locally symmetric spaces of rank $> 1$ depend on the truth of
Schanuel's conjecture (hence are {\it conditional}). Analyzing
the right hand side of equation (\ref{E:3-2}) with the help of the above consequence of
Schanuel's conjecture, we show that if both $\Gamma_1$ and
$\Gamma_2$ can be conjugated into ${\rm SL}_n(\overline{\Q}),$ for
non-weakly commensurable $\gamma_i \in \Gamma_i,$
$\ell_{\Gamma_1}(c_{\gamma_1})$ and $\ell_{\Gamma_2}(c_{\gamma_2})$
are algebraically independent over $\Q.$ Thus, we obtain the
following.

\vskip2mm

\noindent {\bf Proposition.} {\it Let $\Gamma_1$ and $\Gamma_2$ be
discrete torsion-free subgroups of $\cG = G(\R),$ where $G$ is an
absolutely simple $\R$-subgroup of $\mathrm{SL}_n.$ In the case
$\mathrm{rk}_{\R}\: G > 1,$ assume that Schanuel's conjecture holds
and both $\Gamma_1$ and $\Gamma_2$ can be conjugated into
${\rm SL}_n(\overline{\Q}).$ If $\fX_{\Gamma_1}$ and $\fX_{\Gamma_2}$ are
length-commensurable, then the subgroups $\Gamma_1$ and $\Gamma_2$
are weakly commensurable.}

\vskip2mm

We recall that if $\Gamma$ is a lattice in $\cG = G(\R),$ where $G$
is an absolutely simple real algebraic group, not isogenous to
$\mathrm{SL}_2,$ then there exists a real number field $K$ such that $G$
is defined over $K$ and $\Gamma \subset G(K),$ see \cite{Ra},
Proposition 6.6. In particular, if $\mathrm{rk}_{\R}\: G > 1$ and
$\Gamma$ is a lattice in $\cG$ (or, equivalently, $\fX_{\Gamma}$ has
finite volume), then $\Gamma$ can always be conjugated into
${\rm SL}_n(\overline{\Q}),$ so the corresponding assumption in the
proposition is redundant. Theorem A now implies

\vskip2mm

\noindent {\bf Theorem 1.} {\it Assume that $\fX_{\Gamma_1}$ and
$\fX_{\Gamma_2}$ are of finite volume, and let $K_{\Gamma_i}$ denote
the field generated by the traces $\mathrm{Tr}\: \mathrm{Ad}\:
\gamma$ for $\gamma \in \Gamma_i.$ If $\fX_{\Gamma_1}$ and
$\fX_{\Gamma_2}$ are length-commensurable, then $K_{\Gamma_1} =
K_{\Gamma_2}.$}

\vskip2mm

We now turn to arithmetically defined locally symmetric spaces.
Combining Theorems C and D with the above proposition, we obtain the
following.

\vskip2mm

\noindent {\bf Theorem 2.} {\it Each class of length-commensurable
arithmetically defined locally symmetric spaces of $\cG = G(\R)$ is
a union of finitely many commensurability classes. It in fact
consists of a~single commensurability class if $G$ is not of type
$A_n$ $(n > 1),$ $D_{2n+1}$ $(n\geqslant 1),$ or $E_6.$}

\vskip2mm

Next, Theorems $E$ and $F$ imply

\vskip2mm

\noindent {\bf Theorem 3.} {\it Assume that  $\fX_{\Gamma_1}$ and
$\fX_{\Gamma_2}$ are of finite volume, and at least one of them is
arithmetically defined. If they are length-commensurable then
both are arithmetically defined and compactness of one of them
implies the compactness of the other.}

\vskip2mm

We now recall that isospectral compact locally symmetric spaces  have same 
weak length spectrum\:(\cite{PR5}, Theorem 10.1). Combining this fact
with Theorems 2 and 3, we obtain the following results, which
apparently do not follow directly from the spectral theory.

\vskip2mm

\noindent{\bf Theorem 4.} {\it Any two arithmetically defined
compact isospectral locally symmetric spaces of an absolutely simple
real Lie group of type other than $A_n$ $(n
> 1),$ $D_{2n+1}$ $(n\geqslant 1)$, or $E_6$, are commensurable to each other.}

\vskip2mm

\noindent {\bf Theorem 5.} {\it If two compact locally symmetric spaces of an absolutely 
simple Lie group are isospectral, and one of them is
arithmetically defined, then the other is also arithmetically
defined.}

\vskip2mm

Finally, assuming Schanuel's conjecture we can obtain the following
(unpublished) strengthening of the result of \cite{Re} for hyperbolic 2-manifolds.

\vskip2mm

\noindent {\bf Theorem 6.} {\it Let $M_1$ and $M_2$ be
arithmetically defined compact hyperbolic 2-manifolds which are not
commensurable. Let $\mathscr{L}_i$ denote the subfield of $\R$
generated (over $\Q$) by $L(M_i).$ Then $\mathscr{L} :=
\mathscr{L}_1\mathscr{L}_2$ has infinite transcendence degree over
either $\mathscr{L}_1$ or $\mathscr{L}_2.$}

\vskip2mm

It would be interesting to show that a similar statement holds for
arbitrary locally symmetric spaces $\fX_{\Gamma_1}$ and
$\fX_{\Gamma_2}$ assuming that they are not length-commensurable.

\section{Proofs: $p$-adic techniques}

Given two arithmetic subgroups, or, more generally, two Zariski-dense subgroups, the proofs of Theorems A--F ultimately rely on the possibility of
constructing semi-simple elements in one subgroup whose spectra are quite 
different from the spectra of {\it all} semi-simple in the other
subgroup unless certain strong conditions relating these subgroups
hold. These results fit into a broader project of
constructing elements with special properties in a given
Zariski-dense subgroup dealt with in our papers \cite{PR1},
\cite{PR3}-\cite{PR4}. The starting point of this project was the
following question asked independently by G.A.~Margulis and
R.~Spatzier: {\it Let $\Gamma$ be a Zariski-dense arithmetic
subgroup of a simple algebraic group $G.$ Does there exist a regular
semisimple $\gamma \in \Gamma$ such that $\langle \gamma \rangle$ is
Zariski-dense in $T := Z_G(\gamma)^{\circ}?$} It should be pointed
out that the existence of such an element is by no means obvious.
For example, if $\varepsilon \in \C^{\times}$ is any element of
infinite order, then the subgroup $\langle \varepsilon \rangle \times
\langle \varepsilon \rangle  \subset \C^{\times} \times \C^{\times}$
is Zariski-dense, but it contains no Zariski-dense cyclic subgroup.
Elaborating on this observation, one can construct a $\Q$-torus $T$
such that $T(\Z)$ is Zariski-dense in $T,$ but no element of $T(\Z)$
generates a Zariski-dense subgroup of $T.$ Something similar may
also happen in the semi-simple situation. Namely, let $G$ be a
simple $\Q$ group with $\mathrm{rk}_{\R}\: G = 1.$ Then if a
$\Q$-subtorus $T$ of $G$ has a nontrivial decomposition into an
almost direct product  $T = T_1 \cdot T_2$ over $\Q$ (and such
a decomposition exists if $T$ has a nontrivial $\Q$-subtorus),
no element of $T(\Z)$ generates a Zariski-dense subgroup of $T.$ The
latter example shows  that the fact that a given torus contains 
a proper subtorus is an obstruction to the existence of an element with 
the desired property.  So, in \cite{PR1} we singled out tori which were called ``irreducible",  and used them to provide an affirmative answer
to the question of Margulis and Spatzier.

\vskip2mm

\noindent {\bf Definition.} A $K$-torus $T$ is {\it $K$-irreducible}
if it does not contain any proper $K$-subtori.

\vskip2mm

The point is that if $T$ is $K$-irreducible then {\it any} $t \in
T(K)$ of infinite order generates a Zariski-dense subgroup. So,
given a simple group $G$ defined over a number field $K,$ to find a
required element $\gamma$ in a given $S$-arithmetic subgroup
(assuming that the latter is Zariski-dense), it is enough to
construct an irreducible maximal $K$-torus $T$ of $G$ such that the
group $T(\cO_K(S))$ is infinite.

We will now outline a general procedure for constructing irreducible
tori. Let $T$ be a $K$-torus, $\mathscr{G}_T = \Ga(K_T/K)$, where
$K_T$ is the splitting field of $T.$ Then $T$ is $K$-irreducible if
and only if $\mathscr{G}_T$ acts irreducibly on $X(T) \otimes_{\Z}
\Q.$ Now, if $T$ is a maximal $K$-torus of $G,$ then $\mathscr{G}_T$
acts faithfully on the root system $\Phi(G , T),$ which allows us to
identify $\mathscr{G}_T$ with a subgroup of $\mathrm{Aut}(\Phi(G ,
T)).$ If under this identification $\mathscr{G}_T$ contains
the Weyl group $W(G , T)$, then $T$ is $K$-irreducible. Therefore, it would suffice to find a way to construct maximal $K$-tori $T$ of $G$ such
that $\mathscr{G}_T \supset W(G , T).$ For $G = \mathrm{SL}_n,$ one
simply needs to find a polynomial $f(t) \in K[t]$ of degree $n$ with
Galois group $S_n;$ the existence of such polynomials is well-known. Apparently, similar
constructions can be implemented to obtain irreducible $K$-tori in any other group of   
classical type, but an additional difficulty one has to deal with is that it needs to be shown 
that the torus one constructs admits a $K$-embedding into the group. For the groups of exceptional types, an explicit construction appears to be difficult. Our proof of the existence of required tori in
\cite{PR1} was based on so-called ``generic tori".

It was shown by V.E.~Voskresenskii \cite{Vos} that  $G$ has a
maximal torus $\mathscr{T}$ defined over a purely transcendental
extension $\mathscr{K} = K(x_1, \ldots , x_n)$ such that
$\mathscr{G}_{\mathscr{T}} \supset W(G , \mathscr{T}).$ Then, using
Hilbert's Irreducibility Theorem, one can specialize parameters to
get (plenty of) maximal $K$-tori $T$ of $G$ such that
$\mathscr{G}_T \supset W(G , T).$ In fact, we can construct such
tori with prescribed behavior at finitely many places of $K$, using
which it is easy to ensure that for the resulting torus $T$ the
group $T(\cO_K(S))$ is infinite, and then any element $\gamma \in
T(\cO_K(S))$ of infinite order has the desired property. 

\vskip2mm

Some time later, G.A.\,Margulis and G.A.\,Soifer asked us a different
version of the original question which arose in their joint work with
H.\,Abels on the Auslander problem: {\it Let $G$ be a simple real
algebraic group, $\Gamma$ be a finitely generated Zariski-dense  subgroup of 
$G(\R)$. Is there a regular semi-simple element $\gamma$ 
in $\Gamma$ which generates a Zariski-dense subgroup of $T =
Z_G(\gamma)^{\circ}$ and which is also $\R$-regular?} We recall that
$\gamma \in G(\R)$ is $\R$-{\it regular} if the number of
eigenvalues, counted with multiplicity, of modulus 1 of $\mathrm{Ad}
\gamma$ is minimal possible (cf.\,\cite{PRa}). It should be noted
that even the existence of an $\R$-regular element without any
additional requirement in an arbitrary Zariski-dense subgroup
$\Gamma$ is a nontrivial matter: this was established by Benoist and
Labourie \cite{BL} using the multiplicative ergodic theorem, and
then by Prasad \cite{P} by a direct argument; we will not, however,
discuss this aspect here. The real problem is that the above
argument for the existence of a regular semisimle element in
$\Gamma$ which generates a Zariski-dense subgroup of its centralizer
does not extend to the case where $\Gamma$ is not arithmetic. More
precisely, since $\Gamma$ is finitely generated, we can choose a
finitely generated subfield $K$ of $\R$ such that $G$ is
defined over $K$ and $\Gamma \subset G(K).$ Then
we can construct a maximal $K$-torus $T$ of $G$ which is
irreducible over $K.$ However, it is not clear at all why
$T(K)$ should contain an element of $\Gamma$ of infinite order if the latter is not of
``arithmetic type". Nevertheless, the answer to the question of
Margulis and Soifer turned out to be in the affirmative. 

\vskip2mm

\noindent {\bf Theorem 7}\,(\cite{PR2}). {\it Let $G$ be a connected
semi-simple real algebraic group. Then any Zariski-dense
subsemigroup $\Gamma$ of $G(\R)$ contains a regular $\R$-regular
element $\gamma$ such that the cyclic subgroup generated by it is a
Zariski-dense subgroup of the maximal torus $T =
C_G(\gamma)^{\circ}.$}

\vskip2mm

The proof of the theorem, which we will now sketch, used a rather
interesting technique, viz.\,that of  $p$-adic embeddings. We begin by recalling the following proposition. 

\vskip2mm

\noindent {\bf Proposition}\,(\cite{PR3}). {\it Let $\mathcal{K}$ be a finitely
generated field of characteristic zero, $\mathcal{R} \subset
\mathcal{K}$ be a finitely generated ring. Then there exists an
infinite set of primes $\Pi$  such that for each $p \in \Pi$, there
exists an embedding $\varepsilon_p \colon \mathcal{K}
\hookrightarrow \Q_p$ with the property $\varepsilon_p(\mathcal{R})
\subset \Z_p.$}

\vskip2mm

We will only show that $\Gamma$ contains an ``irreducible" element
$\gamma,$ i.e., a regular semi-simple element whose centralizer $T$ is
a $K$-irreducible maximal torus of $G$. For this, we fix a matrix
realization $G \hookrightarrow \mathrm{SL}_n$ and pick a finitely
generated subring $\mathcal{R}$ of $K$ so that $\Gamma \subset G(\cR) :=
G(K)\cap{\rm SL}_n(\mathcal{R}).$ We
then choose a finitely generated field extension $\cK$ of $K$  over which $G$ splits, 
and fix a $\cK$-split maximal torus $T_0$ of $G$. We now let
$C_1, \ldots , C_r$ denote the nontrivial conjugacy classes in the
Weyl group $W(G , T_0).$ Using the above proposition, we pick $r$ primes
$p_1, \ldots , p_r$ such that for each $p_i$ there is an embedding
$\mathcal{K} \hookrightarrow \Q_{p_i}$ for which $\mathcal{R}
\hookrightarrow \Z_{p_i}.$ We then employ results on Galois
cohomology of semi-simple groups over local field to construct, for
each $i = 1, \ldots , r,$ an open set $\Omega_{p_i}(C_i) \subset
G(\Q_{p_i})$ such that any $\omega \in \Omega_{p_i}(C_i)$ is regular
semi-simple and for $T_{\omega} = Z_G(\omega)^{\circ},$ the Galois
group $\mathscr{G}_{T_{\omega}}$ contains an element from the image
of $C_i$ under the natural identification $[W(G , T_0)] \simeq [W(G ,
T_{\omega})]$, where for a maximal torus $T$ of $G$, $[W(G,T)]$ is the set of conjugacy classes in the Weyl group $W(G,T)$.  To conclude the argument, we show that $$\bigcap_{i =
1}^r \left(\Gamma \cap \Omega_{p_i}(C_i)\right) \neq \emptyset,$$
and any element $\gamma$ of this intersection has the property that
for $T = Z_G(\gamma)^{\circ},$ the inclusion $\mathscr{G}_T \supset
W(G , T)$ holds, as required.

\vskip2mm

Theorem 7 was already used in \cite{AMS}. Furthermore, its suitable
generalizations were instrumental in settling a number of questions
about Zariski-dense subgroups of Lie groups posed by Y.~Benoist,
T.J.~Hitchman and R.~Spatzier (cf.\,\cite{PR4}). As we already
mentioned, the elements constructed in Theorem 7 play a crucial role
in the proof of Theorems A--F.

\vskip2mm

We conclude this article with a brief survey of other applications
of $p$-adic embeddings. To our knowledge, Platonov \cite{Pl} was the
first to use $p$-adic embeddings in the context of algebraic groups.
He proved the following.

\vskip2mm

\noindent {\bf Theorem 8} (\cite{Pl}). {\it If $\pi \colon
\widetilde{G} \to G$ is a nontrivial isogeny of connected semi-simple
groups over a finitely generated field $K$ of characteristic zero
then $\pi(\widetilde{G}(K)) \neq G(K).$}

\vskip2mm

It is enough to show that if $\pi \colon \widetilde{T} \to T$ is a
nontrivial isogeny of $K$-tori then $\pi(\widetilde{T}(K)) \neq T(K).$
For this, we pick a finitely generated extension $\cK$ of $K$ so
that $\widetilde{T}$ and $T$ split over $\cK,$ and
every element of $\mathrm{Ker}\:\pi$ is $\cK$-rational.  Then, using the
above proposition, one finds an embedding $\cK \hookrightarrow \Q_p$
for some $p.$ To conclude the argument, one shows that
$\pi(\widetilde{T}(\cK)) = T(\cK)$ would imply $\pi(\widetilde{T}(\Q_p)) =
T(\Q_p),$ which is obviously false.

\vskip2mm

Another application is representation-theoretic rigidity of groups
with bounded generation (cf.\,\cite{R}, and \cite{PlR}, Appendix A.2).
We recall that an abstract group $\Gamma$ has {\it bounded
generation} if there are elements $\gamma_1, \ldots , \gamma_d \in
\Gamma$ such that
$$\Gamma = \langle \gamma_1 \rangle \cdots \langle \gamma_d
\rangle,$$ where $\langle \gamma_i \rangle$ is the cyclic subgroup
generated by $\gamma_i.$

\vskip2mm

\noindent {\bf Theorem 9} (\cite{R}). {\it Let $\Gamma$ be a group
with bounded generation satisfying the following condition

\vskip1.5mm

\noindent {$(*) \ \ \Gamma'/[\Gamma' ,
\Gamma']$ is finite for every subgroup $\Gamma'$ of  $\Gamma$ of
finite index.

\vskip1.5mm

\noindent Then for any $n \geqslant 1,$ there are only finitely many
inequivalent completely reducible representations $\rho \colon
\Gamma \longrightarrow {\rm GL}_n(\C).$}}

\vskip2mm

The proof is based on the following strengthening  of the above 
proposition: given $\cK$ and $\cR$ as above, there exists an 
infinite set of primes $\Pi$ such that for each $p \in \Pi$ there are 
embeddings $\varepsilon_p^{(i)} \colon \cK \to
\Q_p,$ where $i = 1, 2, \ldots,$  such that
$\varepsilon_p^{(i)}(\cR) \subset \Z_p$ for all $i$,  and
$\varepsilon_p^{(i)}(\cR) \cap \varepsilon_p^{(j)}(\cR)$ consists of
algebraic numbers for all $i\neq j$. The usual argument using representation varieties
show that it is enough to show that for any $\rho \colon \Gamma
\longrightarrow {\rm GL}_n(\C),$ the traces $\Tr\,\rho(\gamma)$ are
algebraic numbers, for all $\gamma \in \Gamma.$ For this we pick a
finitely generated subring $\cR$ of $\C$ for which  $\rho(\Gamma)
\subset {\rm GL}_n(\cR),$ and then fix a prime $p$ for which there are
embeddings $\varepsilon_p^{(i)} \colon \cR \to \Z_p$
as above. Let $\rho^{(i)} \colon \Gamma \longrightarrow {\rm GL}_n(\Z_p)$
be the representation obtained by composing $\rho$ with the
embedding ${\rm GL}_n(\cR) \to {\rm GL}_n(\Z_p)$ induced by
$\varepsilon_p^{(i)}.$ One then observes that bounded generation of
$\Gamma$ implies that for any subgroup $\Gamma'$ of $\Gamma$ of finite index,  
the pro-$p$ completion $\Gamma'_p$ of $\Gamma'$ is a
$p$-adic analytic group. Moreover, $(*)$ implies that the
corresponding Lie algebra is a semi-direct product of a semi-simple
algebra and a nilpotent one where the former acts on the latter
without fixed point. Using the fact that a semi-simple algebra has
only finitely many inequivalent representations in any dimension,
one derives that there are $i \neq j$ such that $\Tr\,
\rho^{(i)}(\gamma) = \Tr\,\rho^{(j)}(\gamma)$ for all $\gamma$ in a
suitable subgroup $\Gamma'$ of $\Gamma$ of finite index. Then it
follows from our construction that the traces $\Tr\, \rho(\gamma)$
are algebraic for $\gamma \in \Gamma',$ and consequently all traces
$\Tr\,\rho(\gamma)$ for $\gamma \in \Gamma$ are algebraic.

\vskip2mm

Finally, we would like to mention the following theorem which
provides a far-reaching generalization of the results of \cite{AMH}
and \cite{MMY}.

\vskip2mm

\noindent {\bf Theorem 10}\,(\cite{PR2}). {\it Let $G$ be a connected
reductive group over an infinite field $K.$ Then no noncentral
subnormal subgroup of $G(K)$ can be contained in a finitely
generated subgroup of $G(K).$}

\vskip2mm

(In fact, a similar result is available in the situation where
$G(K)$ is replaced with the group of points over a semi-local
subring of $K.$) To avoid technicalities, let us assume that $G$ is
absolutely simple, and let $N$ be a noncentral normal (rather than
subnormal) subgroup of $G(K).$ Assume that $N$ is contained in a
finitely generated subgroup of $G(K).$ Then, after fixing a matrix
realization $G \subset \mathrm{SL}_n,$ one can pick a finitely
generated subring $\cR$ of $K$ so that $N \subset G(\cR):= G(K)\cap {\rm SL}_n(\cR).$ Let $\cK$
be a finitely generated field that contains $\cR$, and 
such that $G$ is defined and split over $\cK.$ Now, choose an embedding
$\varepsilon_p \colon \cK \hookrightarrow \Q_p$ so that
$\varepsilon_p(\cR) \subset \Z_p,$ and consider the closures $\Delta
= \overline{N}$ and $\cG = \overline{G(K)}.$ Then $\Delta \subset
G(\Z_p),$ hence it is compact, and at the same time it is normal in $\cG.$ On
the other hand, $\cG$ is essentially $G(\Q_p).$ However, $G(\Q_p)$
does not have any noncentral compact normal subgroups (in fact, the subgroup $G(\Q_p)^+$ of $G(\Q_p)$ is a normal subgroup of finite index which does not
contain {\it any} noncentral normal subgroups). A contradiction.

\vskip5mm

\noindent {\it Acknowledgements.} Both authors were supported in
part by the NSF (grants DMS-0653512 and DMS-0502120) and the
Humboldt Foundation.

\vskip5mm

\bibliographystyle{amsplain}

\begin{thebibliography}{100}

\bibitem[1]{AMS} H.~Abels, G.\,A.~Margulis, G.\,A.~Soifer, {\it The Auslander
conjecture for groups leaving a form of signature $(n-2,2)$
invariant.} Probability in mathematics. Israel J. Math. {\bf 148}(2005),
11--21.


\bibitem[2]{AMH} S.~Akbari, M.~Mahdavi-Hezavehi, {\it Normal
subgroups of $GL_n(D)$ are not finitely generated,} Proc. AMS, {\bf
128}(1999), 1627-1632.

\bibitem[3]{Ba} A.~Baker, {\it Transcendental Number Theory}, 2nd
edition, Cambridge Mathematical Library, Cambridge Univ.\,Press,
1990.


\bibitem[4]{BL} Y.~Benoist, F.~Labourie, {\it Sur les
diff\'eomorphismes d'Anosov affines \`a feuilletages stable et
instable diff\'erentiables,} Invent. math. {\bf 111}(1993), 285-308.


\bibitem[5]{BB} L.~B\'erard-Bergery, {\it Laplacien et
g\'eod\'esiques ferm\'ees sur les formes d'espace hyperbolique
compactes,} S\'eminaire Bourbaki, 24\`eme ann\'ee (1971/1972),
Exp.\,No.\,406, 107-122, Lect. Notes in Math. {\bf 317},
Springer-Verlag, 1973.



\bibitem[6]{CHLR} T.~Chinburg, E.~Hamilton, D.\,D.~Long and A.\,W.~Reid,
{\it Geodesics and commensurability classes of arithmetic hyperbolic
3-manifolds,} preprint (July 2006).


\bibitem[7]{Ga} R.~Gangolli, {\it The length spectra of some compact
manifolds}, J.\,Diff\,Geom.\,{\bf 12}(1977), 403-424.



\bibitem[8]{LMNR} C.\,J.~Leninger, D.\,B.~McReynolds, W.~Neumann,
and A.\,W.~Reid, {\it Length and eigenvalue equivalence,} preprint
(December 2006).


\bibitem[9]{LSV} A.~Lubotzky, B.~Samuels and U.~Vishne, {\it
Division algebras and noncommensurable isospectral manifolds,} Duke
Math.\,J.\,{\bf 135}(2006), 361-379.

\bibitem[10]{Mc} H.\,P.~McKean, {\it The Selberg trace formula as
applied to a compact Riemann surface,} Comm.\,Pure Appl.\,Math. {\bf
25}(1972), 225-246.

\bibitem[11]{MMY} M.~Mahdavi-Hezavehi, M.G.~Mahmudi, S.~Yasamin,
{\it Finitely generated subnormal subgroups of $GL_n(D)$ are
central,} J.~Algebra, {\bf 225}(2000), 517-521.

\bibitem[12]{Pl} V.\,P.~Platonov, {\it Dieudonn\'e's conjecture and
the nonsurjectivity on $k$-points of coverings of algebraic groups,}
Soviet Math. Dokl. {\bf 15}(1974), 927-931.

\bibitem[13]{PlR} V.\,P.~Platonov, A.\,S.~Rapinchuk, {\it Algebraic
Groups and Number Theory,} Academic Press, 1994.




\bibitem[14]{P} G.~Prasad, {\it $\mathbf{R}$-regular elements in
Zariski-dense subgroups,} Quart. J. Math. Oxford Ser. (2) {\bf
45}(1994), 541-545.

\bibitem[15]{PRa} G.~Prasad, M.\,S.~Raghunathan, {\it Cartan
subgroups and lattices in semi-simple groups,} Ann. of Math. (2)
{\bf 96}(1972), 296-317.


\bibitem[16]{PR0} G.~Prasad, A.\,S.~Rapinchuk, {\it Computation of
the metaplectic kernel,} Publ. math. IHES, {\bf 84}(1996), 91-187.


\bibitem[17]{PR1} G.~Prasad, A.\,S.~Rapinchuk,
{\it Irreducible tori in semisimple groups,}
Intern.\,Math.\,Res.\,Notices 2001, {\bf 23}, 1229-1242; Erratum,
{\it ibid} 2002, {\bf 17}, 919-921.

\bibitem[18]{PR2} G.~Prasad, A.\,S.~Rapinchuk, {\it Subnormal subgroups of the
groups of rational points of reductive algebraic groups,} Proc. AMS,
130(2002), 2219-2227.

\bibitem[19]{PR3} G.~Prasad, A.\,S.~Rapinchuk, {\it Existence of
irreducible $\R$-regular elements in Zariski-dense subgroups,}
Math.\,Res.\,Letters {\bf 10}(2003), 21-32.

\bibitem[20]{PR4} G.~Prasad, A.\,S.~Rapinchuk, {\it Zariski-dense
subgroups and transcendental number theory,} Math.\,Res.\,Letters
{\bf 12}(2005), 239-249.

\bibitem[21]{PR5} G.~Prasad, A.S.~Rapinchuk,  {\it Weakly commensurable arithmetic groups, lengths
of closed geodesics and isospectral locally symmetric spaces,}
preprint posted on arXiv. 

\bibitem[22]{PR6} G.~Prasad, A.S.~Rapinchuk, {\it Local-global principles for embedding of fields with involution into simple algebras with involution}, preprint posted on arXiv. 


\bibitem[23]{Ra} M.\,S.~Raghunathan, {\it Discrete subgroups of Lie
groups,} Springer, 1972.


\bibitem[24]{R} A.\,S.~Rapinchuk, {\it Representations of groups of
finite width,} Soviet Math. Dokl. {\bf 42}(1991), 816-820.

\bibitem[25]{Re} A.\,W.~Reid, {\it Isospectrality and commensurability
of arithmetic hyperbolic 2- and 3-manifolds,} Duke Math.\,J. {\bf
65}(1992), 215-228.

\bibitem[26]{Su} T.~Sunada, {\it Riemann coverings and isospectral
manifolds,} Ann.\,Math.\,(2) {\bf 121}(1985), 169-186.

\bibitem[27]{Vi} M-F.\,Vign\'eras, {\it Varietes Riemanniennes
Isospectrales et non Isometriques}, Ann.\,Math.\,(2){\bf 112}(1980),
21-32.

\bibitem[28]{Vn} E.\,B.~Vinberg, {\it Rings of definition of dense
subgroups of semisimple linear groups,} Math.\,USSR Izvestija, {\bf 5}(1971),
45-55.

\bibitem[29]{Vos} V.\,E.~Voskresenskii, {\it Algebraic Groups and
their Birational Invariants,} AMS, 1998.

\end{thebibliography}

\end{document}